\documentclass[a4paper,12pt]{article}
\usepackage{amssymb,amsthm,amsmath,amsfonts}
\usepackage{graphicx}
\usepackage{cite}
\usepackage{color}

\newtheorem{theorem}{Theorem}[section]

\newtheorem{corollary}[theorem]{Corollary}
\newtheorem{proposition}[theorem]{Proposition}

\theoremstyle{definition}
\newtheorem{definition}[theorem]{Definition}

\theoremstyle{remark}
\newtheorem{remark}{Remark}

\def\pf{\par\noindent {\em Proof.}~\par\noindent}

\begin{document}

\title{A quaternionic perturbed fractional $\psi-$Fueter operator calculus}
\small{
\author
{Jos\'e Oscar Gonz\'alez-Cervantes$^{(1)}$ and Juan Bory-Reyes$^{(2)\footnote{corresponding author}}$}
\vskip 1truecm
\date{\small $^{(1)}$ Departamento de Matem\'aticas, ESFM-Instituto Polit\'ecnico Nacional. 07338, Ciudad M\'exico, M\'exico\\ Email: jogc200678@gmail.com\\$^{(2)}$ {SEPI, ESIME-Zacatenco-Instituto Polit\'ecnico Nacional. 07338, Ciudad M\'exico, M\'exico}\\Email: juanboryreyes@yahoo.com
}

\maketitle
\begin{abstract}
Quaternionic analysis offers a function theory focused on the concept of $\psi-$hyperholomorphic functions defined as null solutions of the $\psi-$Fueter operator, where $\psi$ is an arbitrary orthogonal base (called structural set) of $\mathbb H^4$. 

The main goal of the present paper is to extend the results given in \cite{BG2}, where a fractional $\psi-$hyperholomorphic function theory was developed. We introduce a quaternionic perturbed fractional $\psi-$Fueter operator calculus, where Stokes and Borel-Pompeiu formulas in this perturbed fractional $\psi-$Fueter setting are presented.
\end{abstract}


\section{Introduction} 
Fractional calculus, a theory allowing integrals and derivatives of arbitrary real or complex order, attracts the attention of scientists and engineers during the last few decades because of numerous applications in diverse fields of science and engineering. The exponential growth of the
number of publications in this area, avoids us having to make extensive list of references. However, in this brief introduction, reader are referred 
to J. Liouville in a series of papers from 1832-1837, where he defined the first outcast of an operator of fractional integration, see for instance \cite{L} and B. Riemann \cite{R} for the construction of the integral-based Riemann-Liouville fractional integral operator (see \cite{CC}), which has been a valuable cornerstone in fractional calculus ever since. A brief history and exposition of the fundamentals of fractional calculus can be found in \cite{Ro}.

Fractional hyperholomorphic function theory is a very recent topic of research, see \cite{CDOP, CCTT, DM, FRV, FKRV, KV, PBBB, V} for more details. In particular, the interest for consi\-de\-ring fractional Laplace and Dirac type operators is devoted in \cite{Ba, Be, FV1, FV2, RRS}.

Nowadays, quaternionic analysis is regarded as a broadly accepted branch of classical analysis offering a successful generalization of complex holomorphic function theory, the most renowned examples are Sudbery's paper\cite{sudbery}. It relies heavily on results on functions defined on domains in $\mathbb R^4$ with values in the skew field of real quaternions $\mathbb H$ associated to a generalized Cauchy-Riemann operator (the so-called $\psi-$Fueter operator) by using a general orthonormal basis in $\mathbb R^4$ (to be named structural set) $\psi$ of $\mathbb H^4$, see, e.g., \cite{No1}. This theory is centered around the concept of $\psi-$hyperholomorphic functions, see \cite{K, MS, S1} and the references given there.

Stokes and Borel-Pompieu associated to a fractional $\psi-$Fueter operator were studied in \cite{BG2}. Now, we present analogous formula to the previous one associated to perturbed fractional $\psi-$Fueter operators that depend on a vector of complex parameters and a quaternion. It is worth noting that a particular case of these perturbed fractional $\psi-$Fueter operators and their associate representation formulas were obtained by a procedure implemented and applied with great effectiveness by the authors in \cite{BG1, G1}.

The structure of the paper reads as follows: Preliminary's section collects some basic facts and definitions on fractional Riemann-Liouville integro-differential operators and a quaternionic analysis associated to a structural set $\psi$, such as the $\psi$-Fueter operator, Stokes and the Borel-Pompieu formulas. Section 3 is devoted to the main results of the paper. Indeed, we study a quaternionic analysis induced by perturbed fractional $\psi-$Fueter operators, where Stokes and Borel-Pompeiu formulas in this perturbed fractional $\psi-$Fueter setting are presented.

\section{Preliminaries} 
Below we give basic definitions and facts on the fractional calculus and quaternionic analysis. These notions will be used throughout the whole paper.
\subsection{Basic definitions and facts on the on Riemann-Liouville fractional integro-differential operators}
One of the most popular (even though it has disadvantages for applications to real world problems) fractional derivatives is the Riemann–Liouville derivative. For completeness, we recall the key definitions and results on Riemann-Liouville fractional integro-differential operators

Given $\alpha\in \mathbb C$ with $\Re \alpha> 0$, let us recall that the Riemann-Liouville integrals of order $\alpha$ of $f\in L^1([a, b], \mathbb R)$, with  $-\infty <a  < b< \infty$, on the left and on the right, are defined by \cite[Definition 2.1]{SKM}. 
$$({\bf I}_{a^+}^{\alpha} f)(x) := \frac{1}{\Gamma(\alpha)} \int_a^x \frac{f(\tau)}{(x-\tau)^{1-\alpha}} d\tau, \quad \textrm{with}  \quad x > a$$
and
$$({\bf I}_{b^-}^{\alpha} f)(x) := \frac{1}{\Gamma(\alpha)} \int_x^b \frac{f(\tau)}{(\tau-x)^{1-\alpha}} d\tau, \quad \textrm{with}  \quad x < b,$$
respectively. 

The space $AC^n([a, b], \mathbb R)$ denotes those functions $f:[a, b]\rightarrow \mathbb R$, which are continuously differentiable on the segment $[a, b]$ up to the order $n-1$ and $f^{(n-1)}$ is supposed to be absolutely continuous on $[a, b]$.

Let $n=[\Re \alpha]+1$, where $[\cdot]$ means the integer part of $\cdot$ and $f\in AC^n([a, b], \mathbb R)$. The fractional derivatives in the Riemann-Liouville sense, on the left and on the right, are given by \cite[Definition 2.2]{SKM}
\begin{align}\label{FracDer} 
(D _{a^+}^{\alpha} f)(x):= \frac{d}{dx^n} \left[ ({\bf I}_{a^+}^{n-\alpha} f)(x)\right]
\end{align}
and
\begin{align} \label{FracDer1}
(D _{b^-}^{\alpha} f)(x):= (-1)^n\frac{d}{dx^n}\left[({\bf I}_{b^-}^{n-\alpha}f)(x)\right] 
\end{align}
respectively. 

It is worth noting that derivatives above exist for $f\in AC^n([a, b], \mathbb R)$. Fractional Riemann-Liouville integral and derivative are linear operators.

Fundamental theorem for Riemann-Liouville fractional calculus \cite{CC} shows that 
\begin{align}\label{FundTheorem}
(D_{a^+}^{\alpha} {\bf I}_{a^+}^{\alpha}f)(x)=f(x) \quad  \textrm{and} \quad (D _{b^-}^{\alpha}  {\bf I}_{b^-}^{\alpha} f)(x) = f(x).
\end{align}

Let us mention an important property of the fractional Riemann-Liouville integral and derivative, see \cite[pag. 1835]{VTRMB}.
\begin{proposition}
\begin{equation}\label{cte}
( D _{a^+}^{\alpha} 1)(x)=\frac{(x-a)^{-\alpha}}{\Gamma[1-\alpha]}  , \quad \forall   x\in [a, b].
\end{equation}
\end{proposition}
\subsection{Rudiments of quaternionic analysis}
Consider the skew field of real quaternions $\mathbb H$ with its basic elements $1, {\bf i},  {\bf j},  {\bf k}$. Thus any element $x$ from $\mathbb H$ is of the form $x=x_0+x_{1} {\bf i}+x_{2} {\bf j}+x_{3} {\bf k}$, $x_{k}\in \mathbb R, k= 0,1,2,3$. The basic elements define arithmetic rules in $\mathbb H$: by definition $ {\bf i}^{2}= {\bf j}^{2}= {\bf k}^{2}=-1$, $ {\bf i}\, {\bf j}=- {\bf j}\, {\bf i}= {\bf k};  {\bf j}\, {\bf k}=- {\bf j}\, {\bf k}= {\bf i}$ and $ {\bf k}\, {\bf i}=- {\bf k}\, {\bf i}= {\bf j}$. For $x\in \mathbb H$ we define the mapping of quaternionic conjugation: $x\rightarrow {\overline x}:=x_0-x_{1}{\bf i}-x_{2}{\bf j}-x_{3}{\bf k}$. In this way it is easy seen that $x\,{\overline x}={\overline x}\,x=x^{2}_{0}+x^{2}_{1}+x^{2}_{2}+x^{2}_{3}$. Note that $\overline {qx}={\overline x}\,\,{\overline q}$ for $q,x\in \mathbb H$. 

The quaternionic scalar product  of $q, x\in\mathbb H$ is given by 
$$\langle q, x\rangle:=\frac{1}{2}(\bar q x + \bar x q) = \frac{1}{2}(q \bar x + x  \bar q).$$

A set of quaternions $\psi=\{\psi_0, \psi_1,\psi_2,\psi_2\}$ is called structural set if $\langle \psi_k, \psi_s\rangle =\delta_{k,s} $, 
for  {$k, s=0,1,2,3$ and} any quaternion $x$  can be rewritten as $x_{\psi} := \sum_{k=0}^3 x_k\psi_k$, where $x_k\in\mathbb R$ for all $k$. Given $q, x\in \mathbb H$ we write
$$\langle q, x \rangle_{\psi}=\sum_{k=0}^3 q_k x_k,$$ 
where $q_k, x_k\in \mathbb R$ for all $k$.

Let $\psi$ an structural set. From now on, we will use the mapping
\begin{equation} \label{mapping}
\sum_{k=0}^3 x_k\psi_k \rightarrow (x_0,x_1,x_2,x_3).
\end{equation}
in essential way.

We have to say something about the set of complex quaternions, which are given by  
$$\mathbb H(\mathbb C(\textsf{i})) =\{q=q_1+  \textsf{i} \ q_2 \ \mid \ q_1,q_2 \in \mathbb H\},$$
where $\textsf{ i}$   is the imaginary unit of $\mathbb C$.
The main difference to the real quaternions is that not all non-zero elements are invertible. There are so-called zero-divisors.

Let us recall that $\mathbb H$ is embedded in $\mathbb H(\mathbb C)$ as follows:
$$\mathbb H =\{q=q_1+   \textsf{i}  \ q_2 \in \mathbb H(\mathbb C)  \ \mid \  q_1,q_2 \in \mathbb H \ \ \textrm{and} \ \ q_2=0\}.$$
The elements of $\mathbb H$ are written in terms of the structural set $\psi$ hence those of $\mathbb H(\mathbb C)$ can be written as $q=\sum_{k=0 }^3 \psi_k q_k,$ where $q_k\in \mathbb C$.
 
Functions $f$ defined in a bounded domain $\Omega\subset\mathbb H\cong \mathbb R^4$ with value in $\mathbb H$ are considered. They may be written as: $f=\sum_{k=0}^3 f_k \psi_k$, where $f_k, k= 0,1,2,3,$ are $\mathbb R$-valued functions in $\Omega$. Properties as continuity, differentiability, integrability and so on, which as ascribed to $f$ have to be posed by all components $f_k$. We will follow standard notation, for example $C^{1}(\Omega, \mathbb H)$ denotes the set of continuously differentiable $\mathbb H$-valued functions defined in $\Omega$.  

The left- and the right-$\psi$-Fueter operators are defined by   
${}^{{\psi}}\mathcal D[f] := \sum_{k=0}^3 \psi_k \partial_k f$ and ${}^{{\psi}}\mathcal  D_r[f] :=  \sum_{k=0}^3 \partial_k f \psi_k$, for all $f \in C^1(\Omega,\mathbb H)$, respectively, where $\partial_k f =\displaystyle \frac{\partial f}{\partial x_k}$ for all $k$. 

Particularly, {if $\partial \Omega$ is a 3-dimensional smooth surface then the Borel-Pompieu formula shows that
\begin{align}\label{BorelHyp}  &  \int_{\partial \Omega}(K_{\psi}(\tau-x)\sigma_{\tau}^{\psi} f(\tau)  +  g(\tau)   \sigma_{\tau}^{\psi} K_{\psi}(\tau-x) ) \nonumber  \\ 
&  - 
\int_{\Omega} (K_{\psi} (y-x) {}^{\psi}\mathcal D [f] (y) + {}^{{\psi}}\mathcal  D_r [g] (y) K_{\psi} (y-x)
     )dy   \nonumber \\
		=  &  \left\{ \begin{array}{ll}  f(x) + g(x) , &  x\in \Omega,  \\ 0 , &  x\in \mathbb H\setminus\overline{\Omega}.                     
\end{array} \right. 
\end{align} 
Differential and integral versions of Stokes' formulas for the $\psi$-hyperholomorphic functions theory are given by \cite{sudbery} 
\begin{align}\label{StokesHyp} d(g\sigma^{{\psi} }_x f) = & \left(g \ {}^{{\psi}}\mathcal  D[f]+ \  {}^{{\psi}}\mathcal D_r[g] f\right)dx,\\
		\int_{\partial \Omega} g\sigma^\psi_x f =  &   \int_{\Omega } \left( g {}^\psi \mathcal  D[f] + {}^{{\psi}}\mathcal  D_r[g] f\right)dx,
\end{align}
for all $f,g \in C^1(\overline{\Omega}, \mathbb H)$}. Here, $d$ stands for the exterior differentiation operator, $dx$ denotes, as usual, the differential form of the 4-dimensional volume in $\mathbb R^4$ and  
$$\sigma^{{\psi} }_{x}:=-sgn\psi \left( \sum_{k=0}^3 (-1)^k \psi_k d\hat{x}_k\right)$$ 
is the quaternionic differential form of the 3-dimensional volume in $\mathbb R^4$ according  to $\psi$, where $d\hat{x}_k  = dx_0 \wedge dx_1\wedge dx_2  \wedge  dx_3 $ omitting factor $dx_k$.  {In addition,} $sgn\psi$ is $1$, or $-1$,  if  $\psi$ and $\psi_{std}:=\{1, {\bf i}, {\bf j}, {\bf k}\}$ have the same orientation, or not, respectively.  

We follow \cite{MS, S1} in assuming that the Cauchy Kernel for the corresponding $\psi$-hyperholomorphic function theory is given by 
 \[ K_{\psi}(\tau- x)=\frac{1}{2\pi^2} \frac{ \overline{\tau_{\psi} - x_{\psi}}}{|\tau_{\psi} - x_{\psi}|^4}.\]
This Cauchy kernel generates the integral operators 
$${}^{\psi}\mathcal T[f](x) = \int_{\Omega} K_{\psi} (y-x) f  (y) dy,$$ 
$${}^{\psi}\mathcal T_r[f](x) = \int_{\Omega}  f  (y)  K_{\psi} (y-x) dy,$$ 
defined for all $f\in L_2(\Omega,\mathbb H)\cup C(\Omega,\mathbb H),$ satisfy 
\begin{align}\label{FueterInv}{}^{\psi}\mathcal D \circ{}^{\psi}\mathcal T[f]= &f , \nonumber \\
{}^{\psi}\mathcal D_r \circ{}^{\psi}\mathcal T_r[f]=&f.  
\end{align} 
We will consider vector parameters $\vec{\alpha} = (\alpha_0, \alpha_1,\alpha_2,\alpha_3) \in\mathbb C^4$ with $ 0 < \Re\alpha_{\ell} < 1$ for $\ell=0,1,2,3$.

Let $a=\sum_{k=0}^3\psi_k a_k, b=\sum_{k=0}^3\psi_k b_k \in \mathbb H$ such that $a_k< b_k$ for all $k$. Write 
\begin{align*}  {J_a^b }:= &  \{  \sum_{k=0}^3\psi_k x_k \in \mathbb H \ \mid \ a_k< x_k < b_k, \  \ k=0,1,2,3\} \\
 = & (a_0,b_0) \times (a_1,b_1) \times (a_2,b_2)  \times (a_3,b_3) ,
\end{align*}
and define $m(J_a^b):=(b_0-a_0) (b_1-a_1)(b_2-a_2)(b_3-a_3)$. 

In what follows, with notation ${J_a^b}$, we assume $a_k < b_k$ for all $k$.

Set {$\vec{\alpha} = (\alpha_0, \alpha_1,\alpha_2,\alpha_3) \in\mathbb C^4$} and $f=\sum_{k=0}^3\psi_k f_k\in AC^1(J_a^b,\mathbb H)$; i.e., $f_k, k= 0,1,2,3,$, the real components of $f$, belongs to $AC^1(\Omega,\mathbb R).$ The mapping $x_j \mapsto f_i(q_0,\dots,x_j,\dots, q_3)$ belongs to $AC^1((a_i, b_i), \mathbb R)$ for each $q\in J_a^b$ and all $i, j=0,1,2,3$.

Now, given $q, x\in  J_a^b $ and $i, j=0,\dots, 3$, the fractional integral of Riemann-Liouville of order { $\alpha_j$} of the mapping 
$x_j \mapsto f_i(q_0,\dots,x_j,\dots, q_3 )$ is defined by: 
$$   ({\bf I}_{a_j^+}^{      {  \alpha_j  }   } f_i)(q_0, \dots, x_j, \dots, q_3)  = \frac{1}{\Gamma(
     {  \alpha_j  }  )} \int_{a_j}^{x_j} \frac{f_i  
(q_0, \dots, \tau_j, \dots, q_3)
}{(x_j-\tau_j)^{1-          {   \alpha_j }  }  } d\tau_j.$$

By the above, as $\displaystyle f=\sum_{i=0}^3 \psi_i f_i$ it follows that 
\begin{align*}  ({\bf I}_{a_j^+}^{    {  \alpha_{j}  } } f )(q_0, \dots, x_j, \dots, q_3) = 
\sum_{i=0}^3 \psi_i ({\bf I}_{{a_j}^+}^{ { \alpha_j }   } f_i)(q_0, \dots, x_j, \dots, q_3)  .
\end{align*}     
for every $f \in AC^1(J_a^b,\mathbb H)$ and $q, x \in J_a^b$.

What is more, the fractional derivative in the Riemann-Liouville sense of the mapping $x_j\mapsto f(q_0, \dots, x_j, \dots, q_3)$  of order $\alpha_j$ is given by
\begin{align*} 
D _{a_j^+}^{      {  \alpha_j  }   } f(q_0, \dots ,  x_j ,\dots ,  q_3) = &
\frac{\partial   }{\partial x_j}
 \sum_{i=0}^3 \psi_i ({\bf I}_{a_j^+}^{      {  \alpha_j  }    } f_i)(q_0, \dots, x_j, \dots, q_3) .
\end{align*}
Note that $({\bf I}_{a_j^+}^{{\alpha_j}}f)$ and $D _{a_j^+}^{{\alpha_j}}f$ are $\mathbb H(\mathbb C)$-valued functions for every $j$. In a similar way we can introduce  $({\bf I}_{b_j^-}^{{\alpha_j}}f)$ and $D _{b_j^-}^{{\alpha_j}}f$. 

Let $f \in AC^1(J_a^b,\mathbb H)$, and $\vec{\alpha} = (\alpha_0, \alpha_1,\alpha_2,\alpha_3) \in\mathbb C^4$ such that $0< \Re \alpha_{\ell}<1$ for $\ell=0,1,2,3$. The fractional $\psi$-Fueter operator of order $\vec{\alpha}$ is defined to be
\begin{align*} 
{}^{\psi}\mathfrak D_a^{\vec{\alpha}}[f] (q,x):= & \sum_{j=0}^3 \psi_j( D _{a_j^+}^{{\alpha_j}}f)(q_0, \dots, x_j , \dots,  q_3)    \\
{}^{\psi}\mathfrak D_{r,a}^{\vec{\alpha}}[f] (q,x):= & \sum_{j=0}^3 ( D _{a_j^+}^{{\alpha_j}}f)(q_0, \dots, x_j , \dots,  q_3) \psi_j    \\
 \end{align*}
for $q, x\in J_a^b$. Note that $q$ is considered a fixed point since the integration and derivation variables are the real components of $x$. Moreover, ${}^{\psi}\mathfrak D_a^{\vec{\alpha}}[f](q,\cdot)$ is a $\mathbb H(\mathbb C)$-valued function.

Particularly, ${}^{\theta}\mathfrak D_a^{\vec{\alpha}}[f] (\xi,q) \mid_{q=\xi}$ can be considered as the fractional derivative of order $\vec \alpha$ at the point $\xi\in J_a^ b$. We shall write ${}^{\theta}\mathfrak D_a^{\vec{\alpha}}[f] (\xi)$ instead of ${}^{\theta}\mathfrak D_a^{\vec{\alpha}}[f] (\xi,q) \mid_{q=\xi}$.

On the other hand, given $f \in AC^1(J_a^b,\mathbb H)$ define 
$${}^{\psi} \mathfrak I_a^{\vec{\alpha}}[f](q,x) = \sum_{j=0}^3\frac{1}{2\Gamma(\alpha_j)} \int_{a_j}^{x_j}\frac{\bar \psi_j 
f(q_0,\dots, \tau_j, \dots, q_3)  + \overline{f(q_0,\dots, \tau_j, \dots, q_3)} \psi_j}{ (x_j- \tau_j)^{1-\alpha_j}} 
d\tau_j  $$  
 and 
	\begin{align*}       {}^{\psi}\mathcal I_a^x [f] (q,x,\vec{\alpha})    
	 =   \int_{J_a^x } \frac{f  
(\tau_0, q_1, \dots, q_3) \frac{(x_0-\tau_0)^{ \alpha_0}}{\Gamma(\alpha_0)}    + \dots+    f  
(q_0, \dots,  q_2, \tau_3)  \frac{(x_3-\tau_3)^{ \alpha_3}  }{\Gamma(\alpha_3)} 
   }{m(J_a^x) }  d\mu_{
\tau},
\end{align*}  
where  $\tau = \sum_{k=0}^3 \psi_k \tau_k$ and $d\mu_{\tau}$ is the differential of volume. 
 
In \cite{BG2} Stokes and Borel-Pompieu type formulas associated to the operators: ${}^{\psi}\mathfrak D_a^{\vec{\alpha}}$ and ${}^{\psi}\mathfrak D_{r,a}^{\vec{\alpha}}$ were considered. The following relation can be found in \cite{BG2} 
\begin{align}\label{fracder}
\displaystyle 
{}^{\psi}\mathfrak D_a^{\vec{\alpha}}[f](q,x) = &  {}^{\psi}\mathcal D_x \circ {}^{\psi}\mathcal I_a^x [f](q, x,\vec{\alpha}) \\
\label{Invfracder}
{}^{\psi}\mathfrak D_a^{\vec{\alpha}} \circ  {}^{\psi} \mathfrak I_a^{\vec{\alpha}}[f](q,x) = & \sum_{j=0}^3 \psi_j f_j(q_0,\dots, x_j, \dots, q_3) \end{align}     
for $f\in AC^{1}({{J_a^b}},\mathbb H)$, where $x= \sum_{j=0}^3 \psi_j x_j$ and $q= \sum_{j=0}^3 \psi_j q_j\in J_a^b$.
 
Particularly, 
\begin{align}\label{equa13}\displaystyle 
{}^{\psi}\mathfrak D_a^{\vec{\alpha}} \circ  {}^{\psi} \mathfrak I_a^{\vec{\alpha}}[f](q,x)\mid_{x=q}= f(q). 
\end{align}

Moreover, $\displaystyle {}^{\bar \psi}\mathcal D_x \circ {}^{\psi}\mathfrak D_a^{\vec{\alpha}}[f](q,x) = \Delta_{\mathbb R^4}\circ {}^{\psi}\mathcal I_a^x [f](q, x,\vec{\alpha})$, where $\Delta_{\mathbb R^4}$ denotes the Laplacian in $\mathbb R^4$ according to the real components of $x$.

If the mapping $x\to \mathcal I_a^{x} [f](q, x,\vec{\alpha})$ belongs to $C^2(J_a^b, \mathbb H)$  for all $q$ and set 
$\vec{\beta} = (\beta_0, \beta_1, \beta_2, \beta_3) \in\mathbb C^4$ with $0< \Re\beta_{\ell} <1$ for $\ell=0,1,2,3$  then we have  
\begin{align}\label{equation01}
 {}^{\psi}\mathfrak D_a^{\vec{\alpha}} \circ  {}^{\psi}\mathfrak D_a^{\vec{\beta}}[f] (q,x) = & \sum_{j=0}^3 \psi_j^2 D _{a_j^+}^{{\alpha_j + \beta_j}} f(q_0, \dots, x_j , \dots,  q_3)   
\end{align} 
and 
\begin{align}\label{equation02} 
{}^{\bar\psi}\mathfrak D_a^{\vec{\alpha}} \circ  {}^{\psi}\mathfrak D_a^{\vec{\beta}}[f](q,x) = & \sum_{j=0}^3 D _{a_j^+}^{{\alpha_j + \beta_j}}f(q_0, \dots, x_j ,\dots, q_3).
\end{align}
Previous identities have a version for the right operators, see \cite{BG2}.

If $h:\Omega\to \mathbb H$ is a fixed function, then  ${}^h M$ and $M^h$ are the operators defined on a set of functions $f$ associated to $\Omega$ by the rule: ${}^h M[ f ]=  hf$  and  $M^h[ f ]= fh$.
\section{Main results}
In this section, we derive the Stokes and Borel-Pompeiu formulas in the perturbed fractional $\psi-$Fueter setting. 
\begin{definition}
Consider $f \in AC^1(J_a^b,\mathbb H)$,   $ \vec{\alpha} = (\alpha_0, \alpha_1,\alpha_2,\alpha_3) \in\mathbb C^4$ with $0< \Re \alpha_{\ell} <1$ for $\ell=0,1,2,3$ and $u,v\in \mathbb H(\mathbb C(\textit{i}))$.   Then denote  
	 	\begin{align*} 
		{}^{\psi}_{u, v}\mathfrak D_{a}^{\vec{\alpha}}[f](q,x) = &  {}^{\psi}\mathfrak D_a^{\vec{\alpha}}[f](q,x) + {}^u M[ f](x)
		 +  {}^v M \circ  {}^{\psi}\mathcal I_a^x [f](q, x, \vec{\alpha})   , \\
		{}^{\psi}\mathfrak D_{r, a, u, v}^{\vec{\alpha}}[f](q,x) = &  {}^{\psi}\mathfrak D_{r,a}^{\vec{\alpha}}[f](q,x) + M^u[f](x) + M^v\circ  {}^{\psi}\mathcal I_a^x [f](q, x, \vec{\alpha}) v   ,\\
{}^{\psi}_u\mathcal H_{a}^x (q,x,\vec{\alpha}) = &  {}^{\psi}\mathcal I_a^x [f] (q,x,\vec{\alpha}) +  {}^{\psi}\mathcal T [ u f ](x) ,\\
{}^{\psi}\mathcal H_{a,u}^x (q,x,\vec{\alpha}) = &  {}^{\psi}\mathcal I_a^x [f] (q,x,\vec{\alpha}) +  {}^{\psi}\mathcal T_r [  f u ](x) .
\end{align*} 
\end{definition}
For to simplify the notations and  computations we only shall work with the basic cases, i.e., with function theory induced by the operators
${}^{\psi}_{u, 0}\mathfrak D_{a}^{\vec{\alpha}}$ and 		${}^{\psi}_{0, v}\mathfrak D_{a}^{\vec{\alpha}}$ with the  respective right-versions of these operators. 
\begin{remark}
Particularly, $	{}^{\psi}_{u,0}\mathfrak D_{a}^{\vec{\alpha}}[f]  (q,x) \mid_{x=q}$ can be considered  as ${}^{\psi}\mathfrak D_a^{\vec{\alpha}}[f](q)  + uf (q)$. Similar conclusion can be obtained for  ${}^{\psi}\mathfrak D_{r ,a, u, 0}^{\vec{\alpha}}[f]$, $   	{}^{\psi}_{0,v}\mathfrak D_a^{\vec{\alpha}}[f]$ and for ${}^{\psi}\mathfrak D_{r,a,0 ,v}^{\vec{\beta}}[f]$.

On the other hand, note that
\begin{align*}     
 {}^{\psi}_u\mathcal H_{a}^x [f](q,x,\vec{\alpha})=&\\   
	     \int_{J_a^x }  \sum_{k=0}\frac{ 2\pi^2 |\tau_{\psi} - x_{\psi}|^4 (x_k-\tau_k)^{ \alpha_k} f(q_{\tau,k}) + \Gamma(\alpha_k)  m(J_a^x) ( \tau_{k} - x_{k})\overline{\psi}_k u f(\tau)  }{  2\pi^2\Gamma(\alpha_k) |\tau_{\psi} - x_{\psi}|^4 m(J_a^x) }   d\mu_{
\tau},
\end{align*} 
\begin{align*}     
 {}^{\psi}\mathcal H_{a,u}^x [f](q,x,\vec{\alpha})=&\\   
	     \int_{J_a^x }  \sum_{k=0}\frac{ 2\pi^2 |\tau_{\psi} - x_{\psi}|^4 (x_k-\tau_k)^{ \alpha_k} f(q_{\tau,k}) + \Gamma(\alpha_k)  m(J_a^x) ( \tau_{k} - x_{k}) f(\tau) u \overline{\psi}_k  }{  2\pi^2\Gamma(\alpha_k) |\tau_{\psi} - x_{\psi}|^4 m(J_a^x)} d\mu_{\tau},
\end{align*} 
where $q_{\tau,0} = (\tau_0, q_1, \dots, q_3)$; $q_{\tau,1} =(q_0, \tau_1, \dots, q_3)$ ...
\end{remark} 

\begin{proposition}\label{propFRACD}
Consider $f \in AC^1(J_a^b,\mathbb H)$, $\vec{\alpha} = (\alpha_0, \alpha_1,\alpha_2,\alpha_3) \in\mathbb C^4$ with $0< \Re \alpha_{\ell} <1$ for $\ell=0,1,2,3$ and $u,v\in \mathbb H(\mathbb C(\textit{i}))$. Then 

\noindent	
1. Operators ${}^{\psi}_{u,0}\mathfrak D_{a}^{\vec{\alpha}}$ and ${}^{\psi}\mathfrak D_{r,a, u, 0}^{\vec{\alpha}}$ satisfy the relations

\noindent			
(a) $\displaystyle {}^{\psi}\mathcal D_x \circ   {}^{\psi}_u\mathcal H_{a}^x [f](q,x,\vec{\alpha})    =  {}^{\psi}_{u,0}\mathfrak D_{a}^{\vec{\alpha}}[f](q,x)$.

\noindent	
(b) $\displaystyle {}^{\psi}\mathcal D_{r,x} \circ   {}^{\psi}\mathcal H_{a,u}^x [f](q,x,\vec{\alpha})    =  {}^{\psi}\mathfrak D_{r,a,u,0}^{\vec{\alpha}}[f](q,x)$.

\noindent	
(c) $$\displaystyle {}^{\psi}_{u,0}\mathfrak D_a^{\vec{\alpha}} \circ  {}^{\psi} \mathfrak I_a^{\vec{\alpha}}[f](q,x) = \sum_{j=0}^3 \psi_j f_j(q_0,\dots, x_j, \dots, q_3)  + u \  {}^{\psi} \mathfrak I_a^{\vec{\alpha}}[f](q,x)$$
and 
$$\displaystyle {}^{\psi}\mathfrak D_{r,a,u,0}^{\vec{\alpha}} \circ {}^{\psi} \mathfrak I_a^{\vec{\alpha}}[f](q,x) = \sum_{j=0}^3 \psi_j f_j(q_0,\dots, x_j, \dots, q_3) + {}^{\psi} \mathfrak I_a^{\vec{\alpha}}[f](q,x)u,$$

\noindent 
where $x= \sum_{j=0}^3 \psi_j x_j$ and $q= \sum_{j=0}^3 \psi_j q_j\in J_a^b$. Particularly, 
$$\displaystyle {}^{\psi}_{u,0}\mathfrak D_a^{\vec{\alpha}} \circ  {}^{\psi} \mathfrak I_a^{\vec{\alpha}}[f](q,x)\mid_{x=q}= f(q) + u \ {}^{\psi} \mathfrak I_a^{\vec{\alpha}}[f](q,q)$$
and 
$$\displaystyle {}^{\psi} \mathfrak D_{r,a,u,0}^{\vec{\alpha}} \circ  {}^{\psi} \mathfrak I_a^{\vec{\alpha}}[f](q,x)\mid_{x=q}= f(q) +  {}^{\psi} \mathfrak I_a^{\vec{\alpha}}[f](q,q) u.$$

\noindent	
(d)
$${}^{\bar \psi}\mathcal D_x \circ	{}^{\psi}_{u,0}\mathfrak D_{a}^{\vec{\alpha}}[f](q,x) = \Delta_{\mathbb R^4}\circ {}^{\psi}\mathcal I_a^x [f](q, x,\vec{\alpha}) + {}^{\bar \psi}\mathcal D_x [uf](x),$$
$${}^{\bar{\psi}} \mathcal  D_{r,x}  \circ	{}^{\psi}\mathfrak D_{r,a, u,0}^{\vec{\alpha}}[f](q,x) = \Delta_{\mathbb R^4}\circ {}^{\psi}\mathcal I_a^x [f](q, x,\vec{\alpha}) + {}^{\bar \psi}\mathcal D_{r,x} [fu](x),$$

\noindent
where $\Delta_{\mathbb R^4}$ denotes the Laplacian in $\mathbb R^4$ according to the real components of $x$.

\noindent
(e) If the mapping $x\to \mathcal I_a^{x} [f](q, x,\vec{\alpha})$ belongs to $C^2(J_a^b, \mathbb H)$  for all $q$ and set 
$\vec{\beta} = (\beta_0, \beta_1, \beta_2, \beta_3) \in\mathbb C^4$ with $0< \Re\beta_{\ell} <1$ for $\ell=0,1,2,3$  then we have  
\begin{align*} 
{}^{\psi}_{u,0}\mathfrak D_a^{\vec{\alpha}} \circ  {}^{\psi}_{v,0}\mathfrak D_a^{\vec{\beta}}[f] (q,x) = & \sum_{j=0}^3 \psi_j^2 D _{a_j^+}^{{\alpha_j + \beta_j}} f(q_0, \dots, x_j , \dots,  q_3) \\ 
 & + u {}^{\psi}\mathfrak D_a^{\vec{\beta}}[f] (q,x) + {}^{\psi}\mathfrak D_a^{\vec{\alpha}} [vf] (q,x)+ uv f(x)   ,
\end{align*} 
\begin{align*} 
{}^{\psi}_{u,0}\mathfrak D_a^{\vec{\alpha}} \circ  {}^{\psi}\mathfrak D_{r,a,v,0}^{\vec{\beta}}[f] (q,x) = & \sum_{j=0}^3 \psi_j D _{a_j^+}^{{\alpha_j + \beta_j}} f(q_0, \dots, x_j , \dots,  q_3) \psi_j \\ 
 & + u {}^{\psi}\mathfrak D_{r,a}^{\vec{\beta}}[f] (q,x) + {}^{\psi}\mathfrak D_a^{\vec{\alpha}} [f] (q,x)v+  u f(x) v  
\end{align*} 
\begin{align*} 
{}^{\bar \psi}_{u,0}\mathfrak D_a^{\vec{\alpha}} \circ  {}^{\psi}_{v,0}\mathfrak D_a^{\vec{\beta}}[f] (q,x) = & {}^{\bar \psi}\mathfrak D_a^{\vec{\alpha}} \circ  {}^{\psi} \mathfrak D_a^{\vec{\beta}}[f] (q,x)   \\ 
 & +     u  {}^{\psi}\mathfrak D_a^{\vec{\beta}}[f] (q,x) + {}^{\bar\psi}\mathfrak D_a^{\vec{\alpha}} [vf] (q,x)+ uv f(x) \\
= & \sum_{j=0}^3 D _{a_j^+}^{{\alpha_j + \beta_j}}f(q_0, \dots, x_j ,\dots, q_3) \\ 
 & +     u  {}^{\psi}\mathfrak D_a^{\vec{\beta}}[f] (q,x) + {}^{\bar\psi}\mathfrak D_a^{\vec{\alpha}} [vf] (q,x)+ uv f(x)   , 
\end{align*} 
\begin{align*} 
{}^{\bar \psi}\mathfrak D_{r,a,u,0}^{\vec{\alpha}} \circ  {}^{\psi}\mathfrak D_{r,a,v,0}^{\vec{\beta}}[f] (q,x) = & 
{}^{\bar \psi}\mathfrak D_{r,a }^{\vec{\alpha}} \circ  {}^{\psi}\mathfrak D_{r,a }^{\vec{\beta}}[f] (q,x)  \\ 
 & +        {}^{ \psi}\mathfrak D_{r,a}^{\vec{\beta}}[f] (q,x) u  + {}^{\bar\psi}\mathfrak D_{r,a}^{\vec{\alpha}} [fv] (q,x)+  f(x)  v u \\
= & \sum_{j=0}^3 D _{a_j^+}^{{\alpha_j + \beta_j}}f(q_0, \dots, x_j ,\dots, q_3) \\ 
 & +        {}^{ \psi}\mathfrak D_{r,a}^{\vec{\beta}}[f] (q,x) u  + {}^{\bar\psi}\mathfrak D_{r,a}^{\vec{\alpha}} [fv] (q,x)+  f(x)  v u.
\end{align*} 
Note that, for $\vec{\alpha}=\vec{\beta}$ and $v=\bar u $ the above formula drawn the fact that the fractional $\psi$-Fueter operator of order $\displaystyle\frac{1+\vec{\alpha}}{2}$ factorizes a fractional $\psi$-Laplace operator defined by ${}^{\psi}\Delta_a^{\vec{\alpha}}:=\sum_{j=0}^3 D _{a_j^+}^{{1+\alpha_j}}$.

\noindent	
2. The operators ${}^{\psi}_{0,u}\mathfrak D_{a}^{\vec{\alpha}}$ and $ {}^{\psi}\mathfrak D_{r,a,0,v}^{\vec{\beta}}$ share the following properties
 
\noindent	
(a) ${}^{\psi}\mathcal D_{x} (e^{\langle u, x \rangle_{\psi} }
 {}^{\psi}\mathcal I_a^x [f](q, x, \vec{\alpha})) = e^{\langle u, x \rangle_{\psi}} \
{}^{\psi}_{0,u}\mathfrak D_{a}^{\vec{\alpha}}[f](q,x)$. 
 
\noindent	
(b) ${}^{\psi}\mathcal D_{r,x} (e^{\langle u, x \rangle_{\psi} }
 {}^{\psi}\mathcal I_a^x [f](q, x, \vec{\beta})  ) =  e^{\langle u, x \rangle_{\psi}} \ {}^{\psi}\mathfrak D_{r,a,0,u} ^{\vec{\beta}}[f](q,x)$.

\noindent	
(c) $ {}^{\bar \psi}\mathcal D_{x}( \ e^{\langle u, x \rangle_{\psi}} \ {}^{\psi}_{0,u}\mathfrak D_{a}^{\vec{\alpha}}[f](q,x) \ ) 
   = \Delta_{x}( \ e^{\langle u, x \rangle_{\psi} }
 {}^{\psi}\mathcal I_a^x [f](q, x, \vec{\alpha}) \ )$. 
 
\noindent	 
(d) ${}^{\bar \psi}\mathcal D_{r,x} ( \ e^{\langle u, x \rangle_{\psi}} \ 
{}^{\psi}\mathfrak D_{r,a,0,u}^{\vec{\beta}}[f](q,x) \ )
   =  \Delta_{x} ( \  e^{\langle u, x \rangle_{\psi} }
 {}^{\psi}\mathcal I_a^x [f](q, x, \vec{\beta}) \ )$.

\noindent	
(e) ${}^{\bar \psi}_{0,u}\mathfrak D_{a}^{\vec{\alpha}}\circ {}^{\psi}_{0,v}\mathfrak D_{a}^{\vec{\beta}}[f] (q,x) = {}^{\bar \psi}{\mathfrak D}_a^{\vec{\alpha}} \circ {}^{\psi}\mathfrak D_a^{\vec{\beta}}[f](q,x) +{}^{\bar \psi}{\mathfrak D}_a^{\vec{\alpha}} \left(v \ {}^{\psi}\mathcal I_a^x [f](q, x, \vec{\beta}) \right) \\
	     +  u\   {}^{\bar \psi}\mathcal I_a^x [ \  {}^{\psi}\mathfrak D_a^{\vec{\beta}}[f]\ ](q,x,\alpha)
+   uv  \ {}^{\bar \psi}\mathcal I_a^x \left[ \ {}^{\psi}\mathcal I_a^x [f]  (q, x, \vec{\beta}  ) \right] (q, x, \vec{\alpha})$.
\end{proposition}
 
\pf

\noindent
Part 1(a)
$$\displaystyle {}^{\psi}\mathcal D_x \circ  \left({}^{\psi}\mathcal I_a^x  [f](q, x,\vec{\alpha}) +  {}^{\psi}\mathcal T [u f](x) \right)  =  {}^{\psi}\mathfrak D_a^{\vec{\alpha}}[f](q,x) + uf(x) =  {}^{\psi}_{u,0}\mathfrak D_{a}^{\vec{\alpha}}[f](q,x).$$
\medskip

\noindent	
Part 1(b) 
$$\displaystyle 
  {}^{\psi}\mathcal D_{r,x} \circ  \left({}^{\psi}\mathcal I_a^x  [f](q, x,\vec{\alpha}) + 
	{}^{\psi}\mathcal T_r [ f u](x) \right)  =  {}^{\psi}\mathfrak D_{r,a}^{\vec{\alpha}}[f](q,x) + f(x)u 
	=  {}^{\psi}\mathfrak D_{r,a,u,0}^{\vec{\alpha}}[f](q,x).$$
\medskip

\noindent
Part 1(c) These facts are consequences of  \eqref{Invfracder}  and  \eqref{equa13}.
\medskip

\noindent	
Part 1(d) From short computations and use of (a) and (b).
\medskip

\noindent
Part 1(e) It is a matter of direct computation using \eqref{equation01} and \eqref{equation02}.
\medskip

\noindent
Part 2(a)
\begin{align*}{}^{\psi}\mathcal D_{x} (e^{\langle u, x \rangle_{\psi} }
 {}^{\psi}\mathcal I_a^x [f](q, x, \vec{\alpha})  )  = & \left( u {}^{\psi}\mathcal I_a^x [f](q, x, \vec{\alpha}) +  
{}^{\psi}\mathfrak D_a^{\vec{\alpha}}[f](q,x) 
\right) e^{\langle b, x \rangle_{\psi}}   \\
=&  
{}^{\psi}_{0,u}\mathfrak D_{a}^{\vec{\alpha}}[f](q,x) 
 e^{\langle u, x \rangle_{\psi}} 
 \end{align*}
\medskip

\noindent
Part 2(b)  
\begin{align*} 
{}^{\psi}\mathcal D_{r,x} (e^{\langle u, x \rangle_{\psi} }
 {}^{\psi}\mathcal I_a^x [f](q, x, \vec{\beta})  ) = & \left({}^{\psi}\mathcal I_a^x [f](q, x, \vec{\beta}) u +  
{}^{\psi}\mathfrak D_{r,a}^{\vec{\beta}}[f](q,x) 
\right) e^{\langle u, x \rangle_{\psi}}  \\
 =&   
{}^{\psi}\mathfrak D_{r,a, 0, u}^{\vec{\beta}}[f](q,x) 
 e^{\langle u, x \rangle_{\psi}}   
\end{align*}
\medskip

\noindent
Part 2(c) Applies $ {}^{\bar \psi}\mathcal D_{x} $ on both sides of 2(a).
\medskip

\noindent
Part 2(d) Similar to above.
\medskip

\noindent
Part 2(e)  
\begin{align*}  
{}^{\bar \psi}_{0,u}\mathfrak D_{a}^{\vec{\alpha}}\circ {}^{\psi}_{0,v}\mathfrak D_{a}^{\vec{\alpha}}[f] (q,x) 
  = & {}^{\bar \psi}_{0,u}{\mathfrak D}_a^{\vec{\alpha}} \circ  {}^{\psi}\mathfrak D_a^{\vec{\alpha}}[f](q,x) 
 + {}^{\bar  \psi}_{0,u}{\mathfrak  D}_a^{\vec{\alpha}} \left( v {}^{\psi}\mathcal I_a^x [f](q, x, \vec{\alpha})  \right) \\
  = &
{}^{\bar  \psi}{\mathfrak  D}_a^{\vec{\alpha}} \circ  {}^{\psi}\mathfrak D_a^{\vec{\alpha}}[f](q,x) +  u {}^{\bar  \psi}\mathcal I_a^x  [{}^{\psi}\mathfrak D_a^{\vec{\alpha}}[f](q,x) ] (q, x, \vec{\alpha}  ) \\
 &  +{}^{\bar \psi}{\mathfrak  D}_a^{\vec{\alpha}} \left( v {}^{\psi}\mathcal I_a^x [f](q, x, \vec{\alpha})  \right) 
+   u {}^{\bar  \psi}\mathcal I_a^x  [ \ v {}^{\psi}\mathcal I_a^x [f] \ ] (q, x, \vec{\alpha}).
\end{align*}
 \hfill$\square$
\begin{remark}\label{Remark1}  
Items 2(a) and 2(b) of previous proposition are structurally analogous to \cite[Formulas 3.1]{BG1} since the exponential function was used to introduce a displaced in terms of a quaternion in both function theories. 
\end{remark}

\begin{proposition} 
If $\vec{\alpha},\vec{\beta} \in\mathbb C^4$ with $0< \Re\alpha_{\ell}, \Re\beta_{\ell}<1$ for $\ell=0,1,2,3$ and  let $f,g \in AC^1(\overline{J_a^b}, \mathbb H)$ consider $q\in J_a^b$ such that the mappings $x\mapsto {}^{\psi}\mathcal I_a^x [f](q,x, \vec{\alpha})$ and $ x\mapsto {}^{\psi} \mathcal I_a^x [g](q,x, \vec{\beta} )$ belong to  $ C^1(\overline{J_a^b}, \mathbb H(\mathbb C))$. For $u,v\in \mathbb H(\mathbb C(\textit{i}))$ the following formulas hold.
 
\noindent 
1. Stokes type integral formula induced by  ${}^{\psi}_{u,0}\mathfrak D_{a}^{\vec{\alpha}}$ and ${}^{\psi}\mathfrak D_{r,a,v, 0}^{\vec{\beta}}$.  
\begin{align*} & \int_{\partial J_a^b} {}^{\psi}\mathcal H_{a,v}^x [g](q,x,\vec{\beta}) \ \sigma^{{\psi} }_x \ {}^{\psi}_u\mathcal H_{a}^x [f](q,x,\vec{\alpha}) 
\\
= &  \int_{J_a^b } [ \  {}^{\psi}\mathcal H_{a,v}^x [g](q,x,\vec{\beta})  \ {}^{\psi}_{u,0}\mathfrak D_a^{\vec{\alpha}}[f](q, x) +
 \  {}^{\psi}\mathfrak D_{r,a,v,0}^{\vec{\beta}}[g](q, x) \ {}^{\psi}_u\mathcal H_{a}^x [f](q,x,\vec{\alpha})\ ]dx.
\end{align*}

\noindent 
2. Stokes type integral formula induced by ${}^{\psi}_{0,u}\mathfrak D_{a}^{\vec{\alpha}}$ and $ {}^{\psi}\mathfrak D_{r,a,0,v}^{\vec{\beta}}$.   
\begin{align*} &   \int_{\partial J_a^b} 
 {}^{\psi}\mathcal I_a^x [g](q, x, \vec{\beta})  
 \ \delta_{x, u+v}^{\psi} \ {}^{\psi}\mathcal I_a^x [f](q, x, \vec{\alpha}) 
 \\ 
=  &  \int_{J_a^b }  ( \   {}^{\psi} \mathcal I_a^x [g](q, x, \vec{\beta})  \   
{}^{\psi}_{0,u}\mathfrak D_{a}^{\vec{\alpha}}[f](q,x) +  
{}^{\psi}\mathfrak D_{r,a,0,v}^{\vec{\beta}}[g ](q,x) 
 \ {}^{\psi}\mathcal I_a^x [f](q, x, \vec{\alpha})\ ) d\mu_{x,u+v},
\end{align*}
where 
$ \delta_{x, u+v}^{\psi} = e^{\langle u+v, x \rangle_{\psi} } \sigma^{{\psi} }_x $ and  
$ d\mu_{x,u+v}  =   e^{\langle u+v, x \rangle_{\psi} } dx $.   
\end{proposition}

\pf

\noindent  
Part 1. Applying \eqref{StokesHyp} in ${}^{\psi}\mathcal H_{a,v}^x [g](q,x,\vec{\beta})$ and $ {}^{\psi}_u\mathcal H_{a}^x [f](q,x,\vec{\alpha})$.
\medskip

\noindent 
Part 2. Using \eqref{StokesHyp} for $e^{\langle v, x \rangle_{\psi} } \ {}^{\psi}\mathcal I_a^x [g](q, x, \vec{\beta}) $ and  $e^{\langle u, x \rangle_{\psi} } \ {}^{\psi}\mathcal I_a^x [f](q, x, \vec{\alpha})$ \hfill$\square$.  
\begin{corollary}  
If $\vec{\alpha},\vec{\beta} \in\mathbb C^4$ with $0< \Re\alpha_{\ell}, \Re\beta_{\ell}<1$ for $\ell=0,1,2,3$ and  let $f,g \in AC^1(\overline{J_a^b}, \mathbb H)$ consider $q\in J_a^b$ such that the mappings $x\mapsto {}^{\psi}\mathcal I_a^x [f](q,x, \vec{\alpha})$ and $ x\mapsto {}^{\psi} \mathcal I_a^x [g](q,x, \vec{\beta} )$ belong to  $ C^1(\overline{J_a^b}, \mathbb H(\mathbb C))$.   Set $u,v\in \mathbb H(\mathbb C(\textit{i}))$.
\begin{enumerate} 
\item (Cauchy type integral theorem induced by  ${}^{\psi}_{u,0}\mathfrak D_{a}^{\vec{\alpha}}$ and 
 ${}^{\psi}\mathfrak D_{r,a,v,0}^{\vec{\beta}}$ ) If 
$${}^{\psi}_{u,0}\mathfrak D_a^{\vec{\alpha}}[f](q, x) = {}^{\psi}\mathfrak D_{r,a,v,0}^{\vec{\beta}}[g](q, x)=0,  \quad \forall x\in J_a^b, $$
 then \begin{align*}     \int_{\partial J_a^b} {}^{\psi}\mathcal H_{a,v}^x [g](q,x,\vec{\beta}) \ 
 \sigma^{{\psi} }_x \ {}^{\psi}_u\mathcal H_{a}^x [f](q,x,\vec{\alpha})
 =  0.
\end{align*}
\item (Cauchy type integral theorem induced by ${}^{\psi}_{0,u}\mathfrak D_{a}^{\vec{\alpha}}$ and ${}^{\psi}\mathfrak D_{r,a,0,v}^{\vec{\beta}}$)  
If    
$${}^{\psi}_{0,u}\mathfrak D_{a}^{\vec{\alpha}}[f](q,x) =
{}^{\psi}\mathfrak D_{r,a,0,v}^{\vec{\beta}}[g ](q,x)  =0, \quad  \forall x\in J_a^b,$$
 then 
\begin{align*}     \int_{\partial J_a^b} 
 {}^{\psi}\mathcal I_a^x [g](q, x, \vec{\beta}) \   
  \delta_{x, u+v}^{\psi} \ {}^{\psi}\mathcal I_a^x [f](q, x, \vec{\alpha}) 
= 0.
\end{align*}
\end{enumerate}
\end{corollary}
\begin{theorem}\label{B-P-F-D}  
Let $\vec{\alpha},\vec{\beta}\in\mathbb C^4$ with $0< \Re\alpha_{\ell}, \Re\beta_{\ell}<1$ for $\ell=0,1,2,3$ and $f,g \in AC^1(\overline{J_a^b}, \mathbb H)$. Consider $q\in J_a^b$ such that the mappings   
 $x\to {}^{\psi}\mathcal I_a^x [f](q,x,\vec{\alpha})$ and $x\to {}^{\psi} \mathcal I_a^x [g](q,x, \vec{\beta})$, for $x\in J_a^b$, belong to $C^1(\overline{J_a^b}, \mathbb H(\mathbb C(\textit{i})))$. For $u,v\in \mathbb H(\mathbb C(\textit{i}))$ the following formulas hold.
		
1. Borel-Pompieu type formula induced by ${}^{\psi}_{u,0}\mathfrak D_{a}^{\vec{\alpha}}$ and ${}^{\psi}\mathfrak D_{r,a,v,0}^{\vec{\beta}}$.
			\begin{align*}  &  
  \int_{\partial J_a^b}  \left( \mathfrak K^{\vec{\alpha}}_{\psi, a} (\tau-x)     \sigma_{\tau}^{\psi} \
	{}^{\psi}_u\mathcal H_{a}^x [f](q,x,\vec{\alpha})	 +  
	{}^{\psi}\mathcal H_{a,v}^x [g](q,x,\vec{\beta})
  \sigma_{\tau}^{\psi}     \mathfrak K^{\vec{\beta}}_{\psi, a} (\tau-x)   \right)
	\\
&	-
\int_{J_a^b}\left( \mathfrak K^{\vec{\alpha}}_{\psi, a} (y-x)  \
  {}^{\psi}_{u,0}\mathfrak D_a^{\vec{\alpha}}[f](q,y)  +
	{}^{\psi}\mathfrak D_{r,a,v,0}^{\vec{\beta}}[g](q,y) \mathfrak K^{\vec{\alpha}}_{\psi, a} (y-x) \right)
	dy  \\
		=  &      \left\{ \begin{array}{ll} 
	\displaystyle \sum_{i=0}^3 \left( f +g\right)(q_0, \dots, x_i, \dots, q_3)   +\mathcal M^{\psi}_a[f,g](q,x,\vec \alpha, \vec\beta)
		, &    x\in	J_a^b ,  \\ 0 , &  x\in \mathbb H\setminus\overline{J_a^b},                    
\end{array} \right. 
\end{align*} 
 where  
\begin{align*}
\mathfrak K^{\vec{\alpha}}_{\psi, a} (y-x) = & \sum_{i=0}^3 D _{a_i^+}^{\alpha_i}   [ K_{\psi}(\tau-x)] \\	
	\mathcal M^{\psi}_a[f,g](q,x,\vec \alpha, \vec\beta)= & 
   \sum_{{{\begin{array}{r}i,j=0 \\ i\neq j \end{array}}}}^3  
\left( \frac{({\bf I}_{a_j^+}^{\alpha_j} f )(q_0, \dots, x_j, \dots, q_3) }{ \Gamma[\alpha_i] (x-a_i)^{\alpha_i}}  
 +\frac{({\bf I}_{a_j^+}^{\beta_j} g )(q_0, \dots, x_j, \dots, q_3) }{ \Gamma[\beta_i] (x-a_i)^{\beta_i}} \right) \\
 & +     \sum_{i=0}^3  D _{a_i^+}^{\alpha_i}  \circ  {}^{\psi}\mathcal T [ u f ] (x)  +    \sum_{i=0}^3  D _{a_i^+}^{\beta_i}   \circ {}^{\psi}\mathcal T_r [ gv] (x)  .  
\end{align*}
The differential operators $D_{a_i^+}^{\alpha_i} $  and $D _{a_i^+}^{\beta_i} $ are given in terms of the real component $x_i$ of $x$, on both sides for $i=0,1,2,3$.

2. Borel-Pompieu-Type formula induced by ${}^{\psi}_{0,u}\mathfrak D_{a}^{\vec{\alpha}}$ and $ {}^{\psi}\mathfrak D_{r,a,0,v}^{\vec{\beta}} $.
 \begin{align*} 
 &  \int_{\partial J_a^b} 
  \left( {\bf K}^{\psi, \vec \alpha}_u (\tau-x) \ \sigma_{\tau}^{\psi} \   
 {}^{\psi}\mathcal I_a^\tau [f](q, \tau, \vec{\alpha})  + {}^{\psi} \mathcal I_a^\tau [g](q, \tau, \vec{\beta})  
  \ \sigma_{\tau}^{\psi}  \    
	{\bf K}^{\psi, \vec \beta}_v (\tau-x)  \right)  \\
	&   -     
\int_{J_a^b}  \left(   {\bf K}^{\psi, \vec \alpha}_u(y-x) 
 \  {}^{\psi}_{0,u}\mathfrak D_{a}^{\vec{\alpha}}[f](q,y)  + {}^{\psi}\mathfrak D_{r,a,0,v}^{\vec{\beta}}[g](q,y) 
 {\bf K}^{\psi, \vec\beta}_v(y-x) \right)
  dy    \\
		=  &      \left\{ \begin{array}{ll}   \sum_{i=0}^3   (f+g) (q_0, \dots, x_i, \dots, q_3) +     N[f](q,x,\vec \alpha)  +   N [g](q,x,\vec \beta)     , &     x\in 
		J_a^b ,  \\ 0 , &  x\in \mathbb H\setminus\overline{J_a^b},                   
\end{array} \right. 
\end{align*} 
where 
 $${\bf K}^{\psi, \vec\alpha}_u(y-x) = \sum_{i=0}^3  D _{a_i^+}^{\alpha_i}  \left[e^{\langle u, y-x \rangle_{\psi}} K_{\psi} (y-x) \right],$$
and 
  $$ N[f ](q,x, \vec{\alpha}) = \displaystyle  \sum_{{{\begin{array}{r}i,j=0 \\ i\neq j \end{array}}}}^3  
\frac{ ({\bf I}_{a_j^+}^{\alpha_j} f  )
(q_0, \dots, x_j, \dots, q_3) }{ \Gamma[\alpha_i] (x-a_i)^{\alpha_i}}.$$
\end{theorem} 

\pf

\noindent 
Part 1. If we consider in formula \eqref{BorelHyp} the functions $({}^{\psi}\mathcal I_a^x [f](q, x, \vec{\alpha}) + 
   {}^{\psi}\mathcal T [ u f ] (x))$ and $({}^{\psi} \mathcal I_a^x [g](q, x, \vec{\beta}) + {}^{\psi}\mathcal T_r [  gv](x) )
$	then we obtain 
 \begin{align*}  &  \int_{\partial J_a^b} (K_{\psi}(\tau-x)\sigma_{\tau}^{\psi}  
( {}^{\psi}\mathcal I_a^\tau [f](q, \tau, \vec{\alpha}) + 
   {}^{\psi}\mathcal T [ u f ] (\tau) )  \\
 & + ( {}^{\psi} \mathcal I_a^{\tau} [g](q, \tau, \vec{\beta}) +
  {}^{\psi}\mathcal T_r [  gv](\tau) )
  \sigma_{\tau}^{\psi} K_{\psi}(\tau-x) ) \nonumber  \\ 
 - &     
\int_{J_a^b} (K_{\psi} (y-x) {}^{\psi}\mathcal D_{y} 
( {}^{\psi}\mathcal I_a^y [f](q, y, \vec{\alpha}) + 
   {}^{\psi}\mathcal T [ u f ] (y) )  \\
 &	+	
{}^{\psi}\mathcal D_{r,y} 
( {}^{\psi} \mathcal I_a^y [g](q, y, \vec{\beta}) +
  {}^{\psi}\mathcal T_r [  gv](y) )
K_{\psi} (y-x))dy   \nonumber \\
		=  &      \left\{ \begin{array}{ll} ( {}^{\psi}\mathcal I_a^x [f](q, x, \vec{\alpha}) + 
   {}^{\psi}\mathcal T [ u f ] (x) )  +  ( {}^{\psi} \mathcal I_a^x [g](q, x, \vec{\beta}) +
  {}^{\psi}\mathcal T_r [  gv](x) )
 , &     x\in 
		J_a^b ,  \\ 0 , &  x\in \mathbb H\setminus\overline{J_a^b}.                   
\end{array} \right. 
\end{align*} 
Suppose $g=0$ on $J_a^b$ and act $\sum_{i=0}^3 D _{a_i^+}^{\alpha_i}$ upon both sides of the above identity, where $D _{a_i^+}^{\alpha_i} $ is given in terms of the real component $x_i$ of $x$. Applying \eqref{equa13} and combining Fubbini's Theorem with Leibniz formula yields 
\begin{align*}  &  
  \int_{\partial J_a^b}  \left[ \sum_{i=0}^3 D _{a_i^+}^{\alpha_i}   K_{\psi}(\tau-x)  \right]  \sigma_{\tau}^{\psi} 
	( {}^{\psi}\mathcal I_a^\tau [f](q, \tau, \vec{\alpha}) + 
   {}^{\psi}\mathcal T [ u f ] (\tau) )
	\\
		  - &
\int_{J_a^b}  \left[ \sum_{i=0}^3  D _{a_i^+}^{\alpha_i}  K_{\psi} (y-x) \right] 
  {}^{\psi}_{u,0}\mathfrak D_a^{\vec{\alpha}}[f](q,y)dy  \\
		=  &      \left\{ \begin{array}{ll} 
	\begin{array}{lll} 
\displaystyle \sum_{i=0}^3 f (q_0, \dots, x_i, \dots, q_3)  +     \sum_{i=0}^3  D _{a_i^+}^{\alpha_i}    {}^{\psi}\mathcal T [ u f ] (x)  \\
      +  \displaystyle  \sum_{{{\begin{array}{r}i,j=0 \\ i\neq j \end{array}}}}^3  
\frac{({\bf I}_{a_j^+}^{\alpha_j} f )(q_0, \dots, x_j, \dots, q_3) }{ \Gamma[\alpha_i] (x-a_i)^{\alpha_i}} ,
	\end{array}
		 &    x\in	J_a^b ,  \\ 0 , &  x\in \mathbb H\setminus\overline{J_a^b}.                    
\end{array} \right. 
\end{align*}  
For $f=0$ and $g\neq 0$, we can do a similar calculation to resemble the previous identity. The proof is completed by summing these two identities. \hfill$\square$ 
\bigskip

\noindent
Part 2. Using formula \eqref{BorelHyp} for $e^{\langle u, x \rangle_{\psi} } \  {}^{\psi}\mathcal I_a^x [f](q, x, \vec{\alpha})$ and $e^{\langle v, x \rangle_{\psi} } \  {}^{\psi} \mathcal I_a^x [g](q, x, \vec{\beta})$ then we obtain
 \begin{align*}  &  \int_{\partial J_a^b} \left( K_{\psi}(\tau-x)\sigma_{\tau}^{\psi}  
\, e^{\langle u, \tau \rangle_{\psi} }  \  {}^{\psi}\mathcal I_a^\tau [f](q, \tau, \vec{\alpha}) + e^{\langle v, \tau \rangle_{\psi} }  \  {}^{\psi} \mathcal I_a^\tau [g](q, \tau, \vec{\beta})  
  \sigma_{\tau}^{\psi} K_{\psi}(\tau-x)  \right)  \nonumber  \\ 
  & -     
\int_{J_a^b}  \left( \  K_{\psi} (y-x) {}^{\psi}\mathcal D_{y} ( \ e^{\langle u, y \rangle_{\psi} }
 \ {}^{\psi}\mathcal I_a^y [f](q, y, \vec{\alpha})  \ )   \right. \\
 &\left. 	+	
{}^{\psi}\mathcal D_{r,y} 
( \  e^{\langle v, y \rangle_{\psi} } \  {}^{\psi} \mathcal I_a^y [g](q, y, \vec{\beta}) \ )
K_{\psi} (y-x)  \ \right )dy   \nonumber \\
		=  &      \left\{ \begin{array}{ll}  e^{\langle u, x \rangle_{\psi} } \ {}^{\psi}\mathcal I_a^x [f](q, x, \vec{\alpha})   +  e^{\langle v, x \rangle_{\psi} } \ {}^{\psi} \mathcal I_a^x [g](q, x, \vec{\beta}) 
 , &     x\in 
		J_a^b ,  \\ 0 , &  x\in \mathbb H\setminus\overline{J_a^b},                   
\end{array} \right. 
\end{align*}  
Then 
\begin{align*} 
 &  \int_{\partial J_a^b} 
\left(e^{\langle u, \tau \rangle_{\psi} } K_{\psi}(\tau-x)\sigma_{\tau}^{\psi}  \ 
 {}^{\psi}\mathcal I_a^\tau [f](q, \tau, \vec{\alpha})    
  +  {}^{\psi} \mathcal I_a^\tau [g](q, \tau, \vec{\beta})  
  \sigma_{\tau}^{\psi}  \   e^{\langle v, \tau \rangle_{\psi} } K_{\psi}(\tau-x)  \right) \nonumber  \\ 
  & -     
\int_{J_a^b} \left( e^{\langle u, y \rangle_{\psi}} K_{\psi} (y-x)
\ {}^{\psi}_{0,u}\mathfrak D_{a}^{\vec{\alpha}}[f](q,y) + 
  {}^{\psi}\mathfrak D_{r,a,0,v}^{\vec{\beta}}[g](q,y) 
 e^{\langle v, y \rangle_{\psi}} \ K_{\psi} (y-x) \right) dy   \nonumber \\
		=  &      \left\{ \begin{array}{ll}  e^{\langle u, x \rangle_{\psi} } \ {}^{\psi}\mathcal I_a^x [f](q, x, \vec{\alpha})   +  e^{\langle v, x \rangle_{\psi} }  \ {}^{\psi} \mathcal I_a^x [g](q, x, \vec{\beta}) 
 , &     x\in 
		J_a^b ,  \\ 0 , &  x\in \mathbb H\setminus\overline{J_a^b}.                   
\end{array} \right. 
\end{align*}  
For $g=0$ we have  
\begin{align*} 
 &  \int_{\partial J_a^b} 
e^{\langle u, \tau-x \rangle_{\psi} } K_{\psi}(\tau-x)\sigma_{\tau}^{\psi}  \   
 {}^{\psi}\mathcal I_a^\tau [f](q, \tau, \vec{\alpha})    -     
\int_{J_a^b}  e^{\langle u, y-x \rangle_{\psi}} K_{\psi} (y-x)
 {}^{\psi}_{0,u}\mathfrak D_{a}^{\vec{\alpha}}[f](q,y) 
  dy    \\
		=  &      \left\{ \begin{array}{ll}   {}^{\psi}\mathcal I_a^x [f](q, x, \vec{\alpha})    , &     x\in 
		J_a^b ,  \\ 0 , &  x\in \mathbb H\setminus\overline{J_a^b}.                   
\end{array} \right. 
\end{align*}   
Let us apply $\sum_{i=0}^3 D _{a_i^+}^{\alpha_i}$ on both sides of above formula, where $D _{a_i^+}^{\alpha_i} $ is given in terms of the real component $x_i$ of $x$. Moreover, combining \eqref{equa13} with Fubini's theorem and Leibniz's formula we obtain    
\begin{align*} 
 &  \int_{\partial J_a^b} 
  \sum_{i=0}^3  D _{a_i^+}^{\alpha_i} \left( e^{\langle u, \tau-x \rangle_{\psi} } K_{\psi}(\tau-x)\right) \sigma_{\tau}^{\psi}  \  
 {}^{\psi}\mathcal I_a^\tau [f](q, \tau, \vec{\alpha})\\
 & -     
\int_{J_a^b}  \sum_{i=0}^3  D _{a_i^+}^{\alpha_i}  \left( e^{\langle u, y-x \rangle_{\psi}} K_{\psi} (y-x) \right)
 {}^{\psi}_{0,u}\mathfrak D_{a}^{\vec{\alpha}}[f](q,y) 
  dy    \\
		=  &      \left\{ \begin{array}{ll} \displaystyle   \sum_{i=0}^3   f (q_0, \dots, x_i, \dots, q_3) +  \mathcal N[f](q,x,\vec \alpha), &  x\in 
		J_a^b ,  \\ 0 , &  x\in \mathbb H\setminus\overline{J_a^b}.                   
\end{array} \right. 
\end{align*} 
Computation for $f=0$ and $g\neq 0$ becomes similar and the proof falls naturally by taking sum of both identities. 
\hfill$\square$
\begin{corollary}\label{C-T-F} 
Let $\vec{\alpha},\vec{\beta}\in\mathbb C^4$ with $0< \Re\alpha_{\ell}, \Re\beta_{\ell}<1$ for $\ell=0,1,2,3$ and $f,g \in AC^1(\overline{J_a^b}, \mathbb H)$. Consider $q\in J_a^b$ such that the mappings   
 $x\to {}^{\psi}\mathcal I_a^x [f](q,x,\vec{\alpha})$ and $x\to {}^{\psi} \mathcal I_a^x [g](q,x, \vec{\beta})$, for $x\in J_a^b$, belong to $C^1(\overline{J_a^b}, \mathbb H(\mathbb C))$. For $u,v\in \mathbb H(\mathbb C(\textit{i}))$ the following formulas hold.

\noindent
1. (Cauchy type formula induced by ${}^{\psi}_{u,0}\mathfrak D_{a}^{\vec{\alpha}}$ and ${}^{\psi}\mathfrak D_{r,a,v,0}^{\vec{\beta}}$) If
$${}^{\psi}_{u,0}\mathfrak D_a^{\vec{\alpha}}[f](q, x) = {}^{\psi}\mathfrak D_{r,a,v,0}^{\vec{\beta}}[g](q, x)=0, \quad \forall x\in J_a^b$$ then
\begin{align*}  &  
  \int_{\partial J_a^b}  \left( \mathfrak K^{\vec{\alpha}}_{\psi, a} (\tau-x)     \sigma_{\tau}^{\psi} \
	{}^{\psi}_u\mathcal H_{a}^x [f](q,x,\vec{\alpha})	 +  
	{}^{\psi}\mathcal H_{a,v}^x [g](q,x,\vec{\beta})
  \sigma_{\tau}^{\psi}     \mathfrak K^{\vec{\beta}}_{\psi, a} (\tau-x)   \right)\\ 
		=  &      \left\{ \begin{array}{ll} 
	\displaystyle \sum_{i=0}^3 \left( f +g\right)(q_0, \dots, x_i, \dots, q_3)   +\mathcal M^{\psi}_a[f,g](q,x,\vec \alpha, \vec\beta)
		, &    x\in	J_a^b ,  \\ 0 , &  x\in \mathbb H\setminus\overline{J_a^b}.                    
\end{array} \right. 
\end{align*} 
Particularly, for $x=q$ we have
 			 \begin{align*}    
   &  \int_{\partial J_a^b}  \left( \mathfrak K^{\vec{\alpha}}_{\psi, a} (\tau-x)     \sigma_{\tau}^{\psi} \
	{}^{\psi}_u\mathcal H_{a}^x [f](q,q,\vec{\alpha})	 +  
	{}^{\psi}\mathcal H_{a,v}^x [g](q,q,\vec{\beta})
  \sigma_{\tau}^{\psi}     \mathfrak K^{\vec{\beta}}_{\psi, a} (\tau-q)   \right)\\
		&=  4\left( f +g\right)(q) +\mathcal M^{\psi}_a[f,g](q,q,\vec \alpha, \vec\beta).
\end{align*}

\noindent
2. (Cauchy type formula induced by ${}^{\psi}_{0,u}\mathfrak D_{a}^{\vec{\alpha}}$ and $ {}^{\psi}\mathfrak D_{r,a,0,v}^{\vec{\beta}} $)
If 
$$ {}^{\psi}_{0,u}\mathfrak D_{a}^{\vec{\alpha}}[f](q,y) = {}^{\psi}\mathfrak D_{r,a,0,v}^{\vec{\beta}}[g](q,y)=0 , \quad  \forall y\in J_a^b $$
then 
 \begin{align*} 
 &  \int_{\partial J_a^b} 
  \left( {\bf K}^{\psi, \vec \alpha}_v (\tau-x)  \  \sigma_{\tau}^{\psi}  \ 
 {}^{\psi}\mathcal I_a^\tau [f](q, \tau, \vec{\alpha})  + {}^{\psi} \mathcal I_a^\tau [g](q, \tau, \vec{\beta})  
  \sigma_{\tau}^{\psi} \  
	{\bf K}^{\psi, \vec \beta}_v (\tau-x)  \right)  \\
			=  &      \left\{ \begin{array}{ll}   \sum_{i=0}^3   (f+g) (q_0, \dots, x_i, \dots, q_3) +  \mathcal N[f](q,x,\vec \alpha)  + \mathcal N[g](q,x,\vec \beta)     , &     x\in 
		J_a^b ,  \\ 0 , &  x\in \mathbb H\setminus\overline{J_a^b},                   
\end{array} \right. 
\end{align*} 
For $x=q$ we get 
 \begin{align*} 
 &  \int_{\partial J_a^b} 
  \left( {\bf K}^{\psi, \vec \alpha}_u (\tau-x)  \sigma_{\tau}^{\psi}  \ 
 {}^{\psi}\mathcal I_a^\tau [f](q, \tau, \vec{\alpha})  + {}^{\psi} \mathcal I_a^\tau [g](q, \tau, \vec{\beta})  
  \sigma_{\tau}^{\psi}  
	\  {\bf K}^{\psi, \vec \beta}_v (\tau-x)  \right)  \\
			=  &    4(f+g)  (q ) +  \mathcal N[f](q,q,\vec \alpha)  + \mathcal N[g](q,q,\vec \beta).   
\end{align*}  
\end{corollary} 
    
\subsection*{Conclusion}
In this work we presented a generalization of Stokes and Borel-Pompeiu formulas of the fractional quaternionic $\psi-$hyperholomorphic function theory developed in \cite{BG2} to the context of a perturbed fractional $\psi-$Fueter setting. The current results can be regarded as a starting point for interesting future works. For instance, the study of the function theory induced by the operator ${}^{ \psi}_{u,v}\mathfrak D_a^{\vec{\alpha}}$, where $u,v\in \mathbb H$. Moreover, it is desirable to obtain reproducing functions, as well as the associated Stokes and Borel-Pompieu type formulas.  

\section*{Declarations}
\subsection*{Funding} Instituto Polit\'ecnico Nacional (grant number SIP20211188) and CONACYT.}
\subsection*{Conflict of interest} The authors declare that they have no conflict of interest regarding the work reported in this paper.
\subsection*{Author contributions} Both authors contributed equally in this paper and typed, read, and approved the final form of the manuscript.
\subsection*{Availability of data and material} Not applicable
\subsection*{Code availability} Not applicable

\small{
}
 
\end{document}